\documentclass[reqno,11pt]{amsart}
\usepackage{amsmath}
\usepackage{amsthm}
\usepackage{amssymb}
\usepackage{amscd}

\newtheorem{thm}{Theorem}[section]
\newtheorem{prop}[thm]{Proposition}
\newtheorem{cor}[thm]{Corollary}
\newtheorem{lem}[thm]{Lemma}

\theoremstyle{definition}
 \newtheorem{df}[thm]{Definition}
 \newtheorem{rem}[thm]{Remark}

\newcommand{\tht}[2]{\vartheta
    \genfrac[]{0pt}{}{0}{#1} (k^{-1} \Omega, #2)}
\newcommand{\im}{\mathrm{Im}\,}
\newcommand{\re}{\mathrm{Re}\,}
\newcommand{\ep}{\varepsilon}
\newcommand{\cp}[1]{\mathbb{CP}^{#1}}
\newcommand{\dr}{\partial}
\newcommand{\lt}{{\rule[-2.3ex]{0pt}{5.6ex}}^t \!\!\!}

\title{Projective embeddings and Lagrangian fibrations 
       of Kummer varieties}
\author{Yuichi Nohara}
\date{}

\address{Graduate School of Mathematics, Nagoya University,
        Chikusa-ku, Nagoya 464-8602, Japan}
\email{m98014i@math.nagoya-u.ac.jp}

\begin{document}
\maketitle

\begin{abstract}
In this paper, we study asymptotic behavior of projective embeddings
of Kummer varieties given by theta functions, and their amoebas.
We prove that a Lagrangian fibration of the Kummer variety can be 
approximated by moment maps of the projective spaces.
\end{abstract}


\section{Introduction}

Let $(X,L)$ be a polarized projective variety defined over $\mathbb{C}$.
Then for sufficiently large integer $k$,
$X$ can be embedded into a projective space
by holomorphic sections of $L^k$:
\[
  \iota_k : X \hookrightarrow \cp{N_k} 
  = \mathbb{P} H^0(X,L^k)^*.
\]
We consider a natural torus action on $\cp{N_k}$.
Then we have a moment map
\[
  \mu_k : \cp{N_k} \longrightarrow \Delta_k \subset
  \mathrm{Lie} (T^{N_k})^*.
\]
Note that $\mu_k$ is a Lagrangian fibration of $\cp{N_k}$
with respect to the Fubini-Study metric $\omega_{\mathrm{FS}}$.
We denote by $B_k$ the image of $X$ under $\mu_k$.
$B_k$ is called a {\it compactified amoeba}.
We study the limit of $B_k$ and the restriction 
$\pi_k : X \to B_k$ of the moment map $\mu_k$ to $X$ 
as $k \to \infty$.

The amoebas heavily depend on the choice of projective embeddings,
or equivalently the choice of basis of $H^0(X,L^k)$.
Of course, there is no natural choice of basis in general.
However, in several cases, Lagrangian fibrations give 
natural basis of $H^0(X,L^k)$.
The relation can be interpreted in terms of
geometric quantization or mirror symmetry.

A typical example is the case of abelian varieties.
Let $A  = \mathbb{C}^n / \Omega \mathbb{Z}^n + \mathbb{Z}^n$ 
be an abelian variety and $L \to A$ a principally polarization.
Then holomorphic sections of $L^k$ are essentially given by the 
{\it theta functions}.
It is well known that there are some natural choices of 
basis of theta functions.
For example,
\[
  \tht {-b}z \,, \quad b \in \frac 1k \mathbb{Z}^n / \mathbb{Z}^n ,
\]
give a basis of $H^0(A,L^k)$, where
\begin{align*}
  \vartheta\genfrac[]{0pt}{}{a}{b}(\Omega,z) &=  
     \sum_{l \in \mathbb Z^n} 
       e \left( \frac 12 {}^t (l+a)\Omega (l+a) + {}^t (l+a)(z+b) \right),\\
  e(t) &=  \exp (2 \pi \sqrt{-1}t)\, .
\end{align*}
In other words, we have the following isomorphism
\[
  H^0(A,L^k) \cong \bigoplus_{b \in \frac 1k \mathbb{Z}^n / \mathbb{Z}^n}
  \mathbb{C} \cdot b \,,
\]
where the right hand side is the vector space formally spanned by 
$k$-torsion points $b \in \frac 1k \mathbb{Z}^n / \mathbb{Z}^n$.
This isomorphism can be given by the Lagrangian fibration
\[
  \pi : A \longrightarrow T^n \,, \quad 
  z = \Omega x + y \longmapsto y \,.
\]
See \cite{W} for the interpretation in terms of geometric quantization,
and \cite{PZ}, \cite{F} for the interpretation in terms of mirror symmetry.
We consider the projective embeddings defined by the above basis:
\[
  \iota_k :  A \hookrightarrow \cp{k^n-1} \,,
  \quad z \mapsto \left( \tht {-b_1}z : \dots :
  \tht {-b_{k^n}}z \right) \,.
\]
It is easy to see that the restriction 
\[
  \pi_k = \mu_k \circ \iota_k : A \to B_k
\]
of the moment map is invariant under the translations
\[
  \Omega x + y \mapsto \Omega (x+a) + y \,,
  \quad a \in \frac 1k \mathbb{Z}^n / \mathbb{Z}^n \,.
\]
Namely, $\pi_k$ looks ``close'' to $\pi$ for large $k$.
In fact, the author \cite{N} proved that $\pi_k$ converges to $\pi$ 
as $k \to \infty$ as maps between compact metric spaces.

In this paper, we consider the case of Kummer varieties
$X = A/(-1)_A$, where $(-1)_A : A \to A$ is the inverse morphism.
In this case, ample line bundles and Lagrangian fibrations can be
induced from those on $A$.
In particular, we have natural basis of holomorphic sections.
We prove that the Lagrangian fibration of $X$ can be also approximated
by the moment maps of projective space.

This paper is organized as follows.
In Section 2, we summarize the asymptotic behavior of theta functions
we use in this paper.
The precise statement of the main theorem is given in Section 3.
Section 4 is devoted to the convergence theorem of K\"ahler metrics
which is a part of the main theorem.
In Section 5, we prove the remaining parts of the main theorem.


\section{Asymptotic behavior of theta functions}

Let
$A = \mathbb{C}^n / \Omega \mathbb{Z}^n + \mathbb{Z}^n$
be an abelian variety 
and consider the principal polarization $L \to A$ defined by
\[
  L = (\mathbb C^n \times \mathbb C )/(\Omega \mathbb Z^n + \mathbb Z^n)\, ,
\]
where the action of $\Omega \mathbb Z^n + \mathbb Z^n$ 
on $\mathbb C^n \times \mathbb C$ is given by
\[
  (z,\zeta) \mapsto (z+\lambda, 
           e^{\pi {}^t \bar{\lambda} (\mathrm{Im}\,\Omega)^{-1}z
           + \frac{\pi}2 {}^t \bar{\lambda} 
           (\mathrm{Im}\,\Omega)^{-1} \lambda} \zeta)
\]
for $\lambda \in \Omega \mathbb Z^n + \mathbb Z^n$.
Then $L$ is symmetric, i.e.
\[
  (-1)_A^* L \cong L \,.
\]

We denote the flat K\"ahler metric in the class $c_1(L)$ by
\[
  \omega_0 = 
  \frac{\sqrt{-1}}2 \sum h_{\alpha  \beta}
  dz^{\alpha} \wedge d \bar{z}^{\beta} \, ,
\]
where we write $(h_{\alpha \beta}) = (\im \Omega)^{-1}$.
Then $h_0 = \exp ( -\pi {}^t z ( \im \Omega)^{-1} \bar{z} )$
gives a Hermitian metric on $L$ such that
$c_1(L,h_0) = \omega_0$.

Let $T^b$ and $T^f$ be two $n$-dimensional tori
$\mathbb{R}^n/\mathbb{Z}^n$
and identify $A \cong T^b \times T^f$ by  
$z=\Omega x + y \leftrightarrow (x,y)$.
Then the natural projection
\[
  \pi : (A,\omega_0) \longrightarrow T^b, \quad
  \Omega x + y \longmapsto y
\]
is a Lagrangian fibration.

For each integer $k$, we denote the subgroups of $k$-torsion points by
\begin{align*}
  T^b_k &\cong T^f_k \cong \frac 1k \mathbb{Z}^n/\mathbb{Z}^n \, ,\\
  A_k &= T^b_k \times T^f_k \, ,
\end{align*}
and we write
$\{b_i\}_{i = 1, \dots k^2} 
= \frac 1k \mathbb{Z}^n/\mathbb{Z}^n = T^b_k$.
Then 
\begin{align*}
  s_{b_i} (z) & = s_i (z) \\
  &= C_{\Omega} k^{- \frac n4}
  \exp \left(\frac {\pi}2 {}^t z (\im \Omega)^{-1}z \right)
  \cdot \tht{-b_i}{kz} \, ,
  \quad b_i \in T^b_k
\end{align*}
give a basis of $H^0(A,L^k)$, where $
C_{\Omega} = 2^{\frac n4}(\det (\mathrm{Im}\,\Omega))^{\frac 14}$.

By looking at the action of Heisenberg groups,
we have the following (see \cite{N}):
\begin{prop}
The basis $\{ s_i \}$ are orthonormal with respect to
the $L^2$-inner product.
\end{prop}

\begin{lem}\label{claim2}
\begin{multline*}
  s_j(z) = Ck^{\frac n4} \exp 
  \left( \frac{\pi k}2 {}^t z (\im \Omega)^{-1} z \right)\\
  \cdot 
  \exp \left( \frac {\pi k}{\sqrt{-1}} {}^t (z-b_j) \Omega^{-1} (z-b_j)
  \right)
  \left( 1 +  \phi \right)
\end{multline*}
with
\[
  |\phi| = O \left( \frac 1{\sqrt{k}} \right)\,,
  \quad |d \phi | = O(1)\,.
\]
\end{lem}

\begin{proof}
From the definition of the theta function, we have
\begin{align*}
  &\tht{-b}{z} \\
  &\quad = \sum_{l \in \mathbb Z^n} 
       e \left( \frac 1{2k} {}^t l \Omega l + {}^t l(z+b) \right)\\
  &\quad = \sum_{l \in \mathbb Z^n} 
       e \left( \frac 12 \lt \left( \frac l{\sqrt{k}}
       + \sqrt{k} \Omega^{-1} (z-b) \right) \Omega
       \left( \frac l{\sqrt{k}} + \sqrt{k} \Omega^{-1} (z-b) \right)
       \right)\\
  &\quad \phantom{= \sum_{l \in \mathbb Z^n} 
       e \left( \frac 12 \lt \left( \frac l{\sqrt{k}}
       + \sqrt{k} \Omega^{-1} (z-b) \right) \right)}
    \times e \left( - \frac k2 {}^t (z-b)\Omega^{-1} (z-b) \right) \,.
\end{align*}
What we need to prove is that
\[
\sum_{l \in \mathbb Z^n}
       e \left( \frac 12 \lt \left( \frac l{\sqrt{k}}
       + \sqrt{k} \Omega^{-1} (z-b) \right) \Omega
       \left( \frac l{\sqrt{k}} + \sqrt{k} \Omega^{-1} (z-b) \right)
       \right)
\]
has the form $Ck^{\frac n2}(1+ \phi)$.
\begin{align*}
 & \sum_{l \in \mathbb Z^n} 
       e \left( \frac 12 \lt \left( \frac l{\sqrt{k}}
       + \sqrt{k} \Omega^{-1} (z-b) \right) \Omega
       \left( \frac l{\sqrt{k}} + \sqrt{k} \Omega^{-1} (z-b) \right)
       \right)\\
 & \quad = k^{\frac n2}
       \sum_{l \in \frac 1{\sqrt{k}} \mathbb{Z}^n + 
       \sqrt{k} \Omega^{-1} (z-b)}
       \exp \left( \pi \sqrt{-1} {}^t l \Omega l \right) 
       \frac 1{\sqrt{k^n}}\\
 & \quad = k^{\frac n2} \left( \int_{\mathbb{R^n}} 
       \exp \left( \pi \sqrt{-1} {}^t l \Omega l \right)  dl 
       + O \left( \frac 1{\sqrt{k}} \right) \right)\\
 & \quad = Ck^{\frac n2}
  \left(1 + O \left( \frac 1{\sqrt{k}} \right) \right)\, .
\end{align*}
Similarly we have
\begin{align*}
  & d \sum_{l \in \mathbb Z^n} 
       e \left( \frac 12 \lt \left( \frac l{\sqrt{k}}
       + \sqrt{k} \Omega^{-1} (z-b) \right) \Omega
       \left( \frac l{\sqrt{k}} + \sqrt{k} \Omega^{-1} (z-b) \right)
       \right) \\
  & \quad = 2 \pi \sqrt{-1} k^{\frac {n+1}2} 
     \sum_{l \in \frac 1{\sqrt{k}} \mathbb{Z}^n + 
       \sqrt{k} \Omega^{-1} (z-b)}
     \frac 1{\sqrt{k^n}} 
     \exp \left( \pi \sqrt{-1} {}^t l \Omega l \right)
     {}^t l dz\\
  & \quad = 2 \pi \sqrt{-1} k^{\frac {n+1}2} 
    \left( \int_{\mathbb{R}^n} 
    \exp \left( \pi \sqrt{-1} {}^t l \Omega l \right)
     ({}^t l dz) dl 
     + O \left( \frac 1{\sqrt{k}} \right) \right)\\
  & \quad = O \left( k^{\frac n2} \right)
\end{align*}
Lemma \ref{claim2} follows from the above estimates.
\end{proof}

\begin{cor}
  There exists constants $C, c >0$ independent of $k$ such that,
  \[
      |s_i(z)|_{h_0}^2
      \le C k^{\frac n2} e^{-ck \cdot d_{T^b}(y,b_i)^2}\, ,
  \]
  for each $z = \Omega x + y \in X$,
  where $ d_{T^b}$ is a distance on $T^b$ .
\end{cor}


\section{Convergence theorem for Kummer varieties}

In this section, we state a convergence theorem of Lagrangian fibration 
for Kummer varieties $X = A/(-1)_A$.

Let $L \to A$ be the line bundle defined in section 2.
Since $L$ is symmetric, there exists a line bundle $M$ on $X$
satisfying
\[
  p^*M \cong L^2 ,
\]
where $p : A \to X$ is the natural projection
(see \cite{B} for the case of Kummer surfaces and \cite{S}
 for higher dimensional case).
From the fact that $p^* : \mathrm{Pic} (X) \to \mathrm{Pic} (A)$
is injective, we have $p^*M^k \cong L^{2k}$.

From the symmetricity condition for $L$, 
the inverse morphism $(-1)_A : A \to A$
lifts to an involution $(-1)_L : L \to L$.
In our situation, the involution $(-1)_L$ is given by
$(z,\zeta) \mapsto (-z, \zeta)$.
We consider the involution on $H^0(A,L^{2k})$ defined by
\[
  H^0(A,L^{2k}) \longrightarrow H^0(A,L^{2k}) \, ,
  \quad s \longmapsto (-1)_{L^{2k}} s (-1)_A \,.
\]
We denote the subspace of invariant sections under the above involution
by $H^0(A,L^{2k})^+$.
By direct computation, we have
\[
  (-1)_{L^{2k}} s_{b_i} (-1)_A = s_{-b_i}
\]
for each $b_i \in T^b_{2k}$.
In particular, $H^0(A,L^k)^+$ is spanned by
$s_{b_i} + s_{-b_i}$.
Note that 
\[
  \dim H^0(A,L^{2k})^+ = \frac{(2k)^n-2^n}2 + 2^n 
  = 2^{n-1}(k^n+1) \, .
\]
It is easy to see that $p^* : H^0(X,M^k) \hookrightarrow H^0(A,L^{2k})$
gives an isomorphism
\[
  H^0(X,M^k) \cong H^0(A,L^k)^+.
\]

Let $\omega$ be the flat orbifold metric in the class $c_1(M)$:
\[
  \omega = \sqrt{-1} \sum h_{\alpha  \beta}
  dz^{\alpha} \wedge d \bar{z}^{\beta} \, .
\]
Then that natural map 
\[
  \pi : X = A/(-1)_A \longrightarrow B = T^b/(-1),
\]
induced from $\pi : A \to T^b$
is a Lagrangian fibration with respect to $\omega$.

For a section $t \in H^0(X,M^k)$ corresponding to
$s \in H^0(A,L^{2k})^+$, its $L^2$-norm is
\[
  \|t \|^2_{L^2(X,\omega)} = \int_X |t|^2 \frac{{\omega}^n}{n!}
  = \frac 12 \int_A |s|^2 \frac{(2 \omega_0)^n}{n!}
  = 2^{n-1} \|s\|_{L^2(A,\omega_0)}\, .
\]
Therefore
\[
  t_i = \begin{cases}
         {\displaystyle \frac 1{\sqrt{2^n}} (s_{b_i} + s_{-b_i})} 
         \, , \quad 
         & \text{if $b_i \in T^b_{2k} \backslash T^b_2$,}\\
         {\displaystyle \frac 1{\sqrt {2^{n-1}}} s_{b_i}} 
         \, , \quad
         & \text{if $b_i \in T^b_2$}
        \end{cases}
\]
give an orthonormal basis of $H^0(X,M^k)$.

Let
\[
  \iota_k : X \longrightarrow \mathbb {CP}^{N_k}, \qquad
  z \longmapsto (t_i(z))
\]
be the projective embedding defined by the above basis.
Let $T^{N_k} \subset SU(N_k +1)$ be the maximal torus which consists of 
diagonal matrices and consider its natural action on $\mathbb {CP}^{N_k}$.
This action is Hamiltonian and its moment map is given by
\[
  \mu_k (Z^0:\cdots :Z^{N_k})= \frac 1{\sum |Z^i|^2}
        \left(|Z^0|^2, \dots , |Z^{N_k}|^2 \right)\, ,
\]
here we identify the dual of the Lie algebra of $T^{N_k}$ with
\[
  \left\{(\xi_0, \dots , \xi_{N_k}) \in {\mathbb R}^{N_k+1} \, 
       \Bigm| \, \sum \xi_i = 1 \right\}.
\]
Denote the image of $\iota_k (X) \subset \mathbb {CP}^{N_k}$ 
under $\mu_k$ by $B_k$ and consider the restriction
of $\mu_k$ to $X$:
\[
  \pi_k := \mu_k \circ \iota_k : X \longrightarrow B_k \, .
\]
We also put $\omega_k = \frac 1k \iota_k^* \omega_{\mathrm FS}$.
We claim that $\pi_k : (X, \omega_k) \to B_k$
converge to the Lagrangian fibration $\pi : (X,\omega) \to B$
in the ``Gromov-Hausdorff topology''.

For that purpose, we need to define distances on $B$ and $B_k$.
We define a metric on $B$ in such a way that
$\pi : (X,\omega) \to B$ is a Riemannian submersion.
To define a distance on $B_k \subset \Delta_k$, it suffices to
define a metric on $\Delta_k$,
where $\Delta_k = \mu_k(\mathbb {CP}^{N_k})$ 
is the moment polytope of $\mathbb {CP}^{N_k}$.
The metric on $\Delta_k$ is also defined in such a way that
\[
  \mu_k : \left(\mathbb{CP}^{N_k}, \frac 1k \omega_{\mathrm{FS}}
  \right)
  \longrightarrow \Delta_k
\]
is a Riemannian submersion in the interior of $\Delta_k$.
This is equivalent to the following definition.
Consider the restriction $\mu_k : \mathbb {RP}^{N_k} \to \Delta_k$
of the moment map to $\mathbb {RP}^{N_k} \subset \mathbb {CP}^{N_k}$.
This is a $2^{N_k}$-sheeted covering which branches on the boundary of 
$\Delta_k$.
By identifying $\Delta_k$ with a sheet of $\mathbb {RP}^{N_k}$, 
the restriction of the normalized Fubini-Study metric
$\frac 1k \omega_{\mathrm {FS}}$ gives a metric on $\Delta_k$.
The distance on $B_k$ is induced from this metric.

\begin{thm}\label{main}
  $\pi_k : (X,\omega_k) \to B_k$ converges to 
  $\pi : (X,\omega) \to B$ as $k \to \infty$ in the following sense:
  \begin{enumerate}
  \renewcommand{\labelenumi}{(\roman{enumi})}
    \item $\omega_k$ converges to $\omega$ in $C^{\infty}$
          on each compact set in $X \backslash \rm{Sing}(X)$.
          Moreover, the sequence $\{ (X,\omega_k) \}$ 
          of compact Riemanian manifolds converges to $(X,\omega)$ 
          with respect to the  Gromov-Hausdorff topology.
    \item $B_k$ converge to $B$ with respect to the 
          Gromov-Hausdorff topology.
    \item The sequence of maps $\pi_k : X \to B_k$ 
          converges to $\pi : X \to B$ as maps between metric spaces.
  \end{enumerate}
\end{thm}

Before the proof, we recall the definition of Gromov-Hausdorff convergence
and convergence of maps.

First we recall the notion of {\it Hausdorff distance}.
Let $Z$ be a metric space and $X,Y \subset Z$ be two subsets.
We denote the $\varepsilon$-neighborhood of $X$ in $Z$ by $B(X,\ep)$.
Then the Hausdorff distance between $X$ and $Y$ is given by
\[
  d^Z_{\mathrm{H}}(X,Y) = \inf \left\{  \ep > 0 \, \left| \,
  X \subset B(Y, \ep),\,\,
  Y \subset B(X, \ep)
   \, \right\}\right. .
\]

For metric spaces $X$ and $Y$, the {\it Gromov-Hausdorff distance} is 
defined by
\[
  d_{\mathrm{GH}}(X,Y) = \inf \{ d^Z_{\mathrm{H}}(X,Y) \,\, 
                 | \,\, X, Y \hookrightarrow Z 
                   \,\,\text{are isometric embeddings.}\},
\] 
or equivalently,
\[
  d_{\mathrm{GH}}(X,Y) = \inf \bigl\{ d^{X \coprod Y}_{\mathrm{H}}(X,Y) \bigr\},
\]
where the infimum is taken over all distances on the disjoint union
$X \coprod Y$ compatible with those on $X$ and $Y$.

Next we recall the notion of convergence of maps (see also \cite{P}).
Let $f_k : X_k \to Y_k$, $f : X \to Y$ be maps between metric spaces.
Suppose that $X_k$ and $Y_k$ converge to $X$ and $Y$ respectively
with respect to the Gromov-Hausdorff distance.
From the definition of Gromov-Hausdorff distances, 
there exist isometric embeddings
\[
  X, \,  X_k \hookrightarrow Z \, 
  \left(=  X \coprod (\coprod_k X_k) \right), \quad
  Y, \,  Y_k \hookrightarrow W \, 
  \left(=  Y \coprod (\coprod_i Y_k) \right)
\]
such that $X_i$ (resp. $Y_k$) converge to $X$ (resp. $Y$) with 
respect to the Hausdorff distance in $Z$ (resp. W).
In this case, we say that $\{f_i\}$ converges to $f$ if 
for every sequence $x_k \in X_k$ converging to $x \in X$,
$f_k(x_k)$ converges to $f(x)$ in $W$. 


\section{Convergence of K\"ahler metrics}

In this section, we prove the first statement of the main theorem.
In the case of abelian variety, it follows from the theorem of
Tian \cite{T} and Zelditch \cite{Zel}. 
To discuss the case of Kummer varieties, 
we need an orbifold version of Zelditch's theorem.

\begin{thm}[Song \cite{S}]\label{orb-TZ}
Let $(X,\omega)$ be a compact K\"ahler orbifold of dimension $n \ge 2$
with only finite isolated singularities 
${\rm{Sing}} (X) = \{e_j\}_{j = 1}^m$,
and $(M,h) \to X$ be an orbifold Hermitian line bundle with
$c_1(M,h) = \omega$.
For each $k \gg 1$, we take an orthonormal basis 
of $H^0(X,M^k)$ and consider the projective embedding
$\iota_k : X \to \mathbb {CP}^{N_k}$ defined by them.
We put $\omega_k = \frac 1k \iota^*_k \omega_{\mathrm{FS}}$ as before.
Then 
\begin{equation}
  \left\| \omega - \omega_k \right\|_{C^q,z}
  \le C_q \left( \frac 1k + k^{\frac q2} e^{-k \delta r(z)^2}
  \right) \,,
  \label{metric-est}
\end{equation}
where $ \| \cdot \|_{C^q,z}$ denote the $C^q$-norm at $z \in X$ and
$r(z)$ is the distance between $z$ and the singular set.
In particular, $\omega_k$ converges to $\omega$ in $C^{\infty}$
on each compact set in $X \backslash \rm{Sing}(X)$.
\end{thm}

\begin{rem}
 Dai-Liu-Ma \cite{DLM} also proved the similar theorem.
\end{rem}

\begin{thm}\label{GH_for_X}
Under the same assumption as in Theorem \ref{orb-TZ},
$(X,\omega_k)$ converges to $(X,\omega)$ with respect to the Gromov-Hausdorff
topology.
\end{thm}

\paragraph{Proof.}
Let $d$ and $d_k$ be the distances on $X$ defined by $\omega$ and $\omega_k$
respectively.
For each $e \in \mathrm{Sing} (X)$, we put
\[
  D_k(e) := \biggl\{ z \in X \, \biggm| \,
  d(z,e) < \sqrt{\frac {\log k}{\delta k}} \,
  \biggr\} 
\]
and
\[
  U_k :=
  X \backslash \bigcup_{e \in \mathrm{Sing} (X)} D_k(e) \,.
\]
Then 
\begin{equation}
  d_{\mathrm{GH}} \bigl( (X,\omega) , (U_k, \omega) \bigr)
  \le O \left( \sqrt{\frac {\log k}{k}} \right)
  \label{GH-1}
\end{equation}
by definition.

From (\ref{metric-est}), we have
\[
  \left\| \omega - \omega_k \right\|_{C^1} 
  \le \frac Ck 
\]
on $U_k$.
This implies that 
\begin{equation}
  d_{\mathrm{GH}} \bigl( (U_k,\omega) , (U_k, \omega_k) \bigr)
  \le O \left( \frac 1k \right) \,.
  \label{GH-2}
\end{equation}
 
Take $e \in \mathrm{Sing} (X)$ and $z_0 \in X$ 
which is close to $e$.
Let
\[
  \gamma : [0,l] \longrightarrow X
\]
be a geodesic from $z_0$ to $e$
with $| \dot{\gamma} | = 1$. 
Then $d(z_0,e) = l$.
From (\ref{metric-est}), 
\[
  | \dot{\gamma} |_{\omega_k}
  \le \left( 1 + \frac Ck + C e^{- \delta k r(\gamma (t))^2} \right)
  = \left( 1 + \frac Ck + C e^{- \delta k (l-t)^2} \right) \,.
\]
Hence we have
\begin{align*}
  d_k(z_0,e) &\le \int_0^l | \dot{\gamma} |_{\omega_k} dt\\
  &\le \int_0^l
       \left( 1 + \frac Ck + C e^{- \delta k (l-t)^2} \right) dt\\
  &\le \left( 1 + \frac Ck \right) l + \frac C{\sqrt{k}} \,.
\end{align*}
Therefore the diameter of $D_k(e)$ with respect to $\omega_k$
can be bounded by $O \left( \sqrt{\frac {\log k}{k}} \right)$.
In particular, we have
\begin{equation}
  d_{\mathrm{GH}} \bigl( (X,\omega_k) , (U_k, \omega_k) \bigr)
  \le O \left( \sqrt{\frac {\log k}{k}} \right) \,.
  \label{GH-3}
\end{equation}

By combining (\ref{GH-1}), (\ref{GH-2}) and (\ref{GH-3}), we obtain
\[
  d_{\mathrm{GH}} \bigl( (X,\omega) , (X, \omega_k) \bigr)
  \le O \left( \sqrt{\frac {\log k}{k}} \right) \,.
\]
\qed


\section{proof of the main theorem}

\subsection{Proof of Theorem \ref{main} (ii).}
To prove the Gromov-Hausdorff convergence, it suffices to construct
{\it $\ep$-Hausdorff approximations} $\varphi_k : B \to B_k$
for large $k$ (see \cite{F2}).

\begin{df}
  Let $(X,d_X)$, $(Y,d_Y)$ be two compact metric spaces.
  A map $\varphi : X \to Y$ is said to be an 
  $\ep$-Hausdorff approximation
  if the following two conditions are satisfied.
  \begin{enumerate}
    \item The $\ep$-neighborhood of $\varphi (X)$ coincides with $Y$.
    \item For each $x,y \in X$,
      \[
        |d_X(x,y) - d_Y(\varphi(x), \varphi(y)) | < \ep \, .
      \]
  \end{enumerate}
\end{df}

Take $0 \in T^f_2$ and identify $B$ with the ``zero section''
\[
  (\{0 \} \times T^b)/(-1)_A \subset X.
\]
Note that the inclusion $B \hookrightarrow X$ 
is not necessarily isometric.
We put
\[
  \varphi_k := \pi_k|_B : B \longrightarrow B_k.
\]
We prove that $\varphi_k$ is a $C \sqrt{\frac{\log k}k}$-Hausdorff 
approximation,
here $C>0$ is a constant independent of $k$.

For each $b \in {\rm{Sing}}(B) = T^b_2/(-1)$, we denote the
$\sqrt{\frac {\log k}{\delta k}}$-neighborhood of the singular fiber
$\pi^{-1}(b)$ by
\[
  N_{b,k} = \biggl\{ z \in X \, \biggm| \, 
    d(z,\pi^{-1}(b)) \le \sqrt{\frac {\log k}{\delta k}} \biggr\}
\]
and set
\begin{align*}
  X(k) &= X \backslash \bigcup_{b \in \rm{Sing}(B)} N_{b,k} \,,\\
  B(k) &= \pi (X(k))\,,\\
  B_k(k) &= \pi_k (X(k)) \,.
\end{align*}
Then
\[
  d_{\mathrm{H}}^B (B, B(k)) \le 
  O \left( \sqrt{\frac {\log k}{k}} \right)
\]
by definition.

For each $\xi \in T_p \cp{N_k}$, 
we decompose it into vertical and horizontal parts:
\[
  \begin{matrix}
    T_p \cp{N_k} & = & T_{\cp{N_k}/\Delta_k , p} & \oplus &
    (T_{\cp{N_k}/\Delta_k , p})^{\perp}\\
    \xi & = & \xi^V & + & \xi^H 
  \end{matrix}
\]
where $T_{\cp{N_k}/\Delta_k ,p} = \ker d \mu_k$ is the tangent space
of the fiber of $\mu_k : \cp{N_k} \to \Delta_k$ 
and $(T_{\cp{N_k}/\Delta_k , p})^{\perp}$
is its orthogonal complement with respect to the Fubini-Study metric.
Let $(Z^0: \dots : Z^{N_k})$ be the homogeneous coordinate and write
\[
  \log \frac{Z^i}{Z^0} = u^i + \sqrt{-1} v^i \, .
\]
Then $T_{\cp{N_k}/\Delta_k}$ and $(T_{\cp{N_k}/\Delta_k})^{\perp}$
are spanned by $\frac{\dr}{\dr v^i}$'s and $\frac{\dr}{\dr u^i}$'s
respectively.

Let $\gamma : [0,l] \to B_k$ be a curve and take a lift
$\tilde{\gamma} : [0,l] \to X$.
Then the length of $\gamma$ is given by
\[
  \int_0^l \left| \left( \frac d{dt} \tilde{\gamma} 
  \right)^H \right|_{\omega_k} dt
  \, .
\]

We also decompose $TX$ into the vertical and horizontal spaces:
\[
  T_zX = T_{X/B,z} \oplus (T_{X/B,z})^{\perp}\,.
\]
Then the length of a curve $\gamma : [0,l] \to B$ is given by
\[
  \int_0^l \left| \frac d{dt} \tilde{\gamma} 
  \right|_{\omega} dt \, ,
\]
where $\tilde{\gamma} : [0,l] \to X$ is a horizontal lift.
Therefor, to prove the theorem, 
we need to compare the these two decompositions.

\begin{lem}\label{claim1}
Let $z \in X(k)$.
\begin{enumerate}
 \item If $\xi \in T_{X/B,z}$, then
  \[
    \Bigl| d \iota_k (\xi)^H \Bigr| \le \frac C{\sqrt k} |\xi|\,.
  \]
 \item If $\xi \in (T_{X/B,z})^{\perp}$,
  \[
    \Bigl| d \iota_k (\xi)^V \Bigr| \le \frac C{\sqrt k}  |\xi|\,.
  \]
\end{enumerate}
\end{lem}

\begin{proof}
Since
\[
  d(z, \pi^{-1}(b_j)) \ge \sqrt {\frac {\log k}{\delta k}} \quad
  \text{or} \quad 
  d(z, \pi^{-1}(b_j)) \ge \sqrt{\frac {\log k}{\delta k}} 
\]
on $p^{-1}(X(k)) \subset A$, 
\[
  \bigl| s_{b_j}(z) \bigr|_{C^r} \le k^{\frac n4 + r} 
  O \left( \frac 1{\sqrt{k}} \right) \quad \text{or} \quad
  \bigl| s_{-b_j}(z) \bigr|_{C^r} \le k^{\frac n4 + r} 
  O \left( \frac 1{\sqrt{k}} \right)
\]
for $r = 0,1$.
From this and Lemma \ref{claim2}, we have
\begin{multline*}
  t_j(z) = Ck^{\frac n4} \exp 
  \left( \frac{\pi k}2 {}^t z (\im \Omega)^{-1} z \right)\\
  \cdot
  \exp \left( \frac {\pi k}{\sqrt{-1}} {}^t (z-b_j) \Omega^{-1} (z-b_j)
  \right) \left( 1 +  \phi \right) \,.
\end{multline*}
Hence
\begin{align}
  \frac{Z^j}{Z^0} &= \frac{t_j(z)}{t_0(z)} \notag \\
  &= C_j \exp \Bigl( 2 \pi \sqrt{-1} {}^t (b_j - b_0) \Omega^{-1} z
     \Bigr) + O \left( \frac 1{\sqrt{k}} \right) \notag \\
  &= C_j \exp \Bigl( 2 \pi k \sqrt{-1} {}^t (b_j - b_0) 
     (x + \re (\Omega^{-1}) y)
     - 2 \pi k {}^t (b_j - b_0) \im (\Omega^{-1}) y
     \Bigr) \notag \\
  & \phantom{= C_j \exp \Bigl( 2 \pi k \sqrt{-1} {}^t (b_j - b_0) 
     (x + \re (\Omega^{-1}) y)
     - 2 \pi k {}^t (b_j - b_0)\Bigr)==}
  + \phi(z)
  \label{non-homog}
\end{align}
for some constant $C_j$.

Recall that $\frac{\dr}{\dr x^i}$ tangents to fibers.
Thus
$(T_{X/X^-})^{\perp}$ is spanned by $J \frac{\dr}{\dr x^i}$'s,
where $J$ is the complex structure on $X$.
By direct computation, we have
\begin{multline*}
  \left( J \frac{\dr}{\dr x^1} ,\dots , J \frac{\dr}{\dr x^n} \right)
  = \left(\frac{\dr}{\dr x^1} ,\dots ,\frac{\dr}{\dr x^n} \right)
   \Bigl( - \re (\Omega^{-1}) (\im \Omega^{-1})^{-1} \Bigr)\\
   + \left(\frac{\dr}{\dr y^1} ,\dots ,\frac{\dr}{\dr y^n} \right)
   (\im \Omega^{-1})^{-1} 
\end{multline*}
and
\begin{align*}
  \frac{\dr}{\dr x^i} 
    \Bigl( {}^t (b_j - b_0) \im (\Omega^{-1}) y \Bigr) &= 0,\\
  J \frac{\dr}{\dr x^i}
    \Bigl( {}^t (b_j - b_0)(x + \re (\Omega^{-1}) y) \Bigr) &= 0\,.
\end{align*}
This means that
\begin{align*}
  \left| \left( \frac{\dr}{\dr x^i} \right)^H \right|
  &\le \frac C{\sqrt k}
   \left| \left( \frac{\dr}{\dr x^i} \right)^V \right|\, ,\\
  \left| \left( J \frac{\dr}{\dr x^i} \right)^V \right|
  &\le \frac C{\sqrt k}
   \left| \left( J \frac{\dr}{\dr x^i} \right)^H \right|\,.
\end{align*}
\end{proof}

\begin{lem}\label{claim3}
There exist a constant $C>0$ such that the 
$C \sqrt{\frac {\log k}k}$-neighborhood
of $\varphi_k (B)$ coincides with $B_k$.
\end{lem}

\begin{proof}
Let $\gamma : [0,l] \to X(k)$ be a geodesic along a fiber of $\pi$.
Then, from Lemma \ref{claim1},
\begin{align*}
  d_{B_k}(\pi_k(\gamma (0)), \pi_k( \gamma (l)) ) 
  &\le \int_0^l \bigl| d \iota_k (\dot{\gamma})^H \bigr|_{\omega_k} dt\\
  &\le \frac C{\sqrt k} \int_0^l \bigl| \dot{\gamma} \bigr|_{\omega} dt\\
  &\le \frac C{\sqrt k} \,.
\end{align*}
This implies that a $\frac C{\sqrt k}$-neighborhood
of $\varphi_k (B(k))$ contains $B_k(k)$.

On the other hand, the diameter of $\pi_k( N_{b,k})$ 
is bounded by
\begin{align*}
  &\text{diameter of $D_k(b)$} + \text{diameter of $\pi_k (\pi^{-1}(b'))$}\\
  &\quad \le O \left( \frac 1{\sqrt k} \right) + 
   O \left( \sqrt{\frac {\log k}k} \right)
  \le O \left( \sqrt{\frac {\log k}k} \right) \,,
\end{align*}
where $b' \in \partial D_k(b)$.
This means that
\[
  d_{\mathrm{H}}^{B_k}(B_k(k), B_k) \le 
  O \left( \sqrt{\frac {\log k}k} \right) \,.
\]
Hence the $C \sqrt{ \frac {\log k}k}$-neighborhood
of $\varphi_k (B)$ coincides with $B_k$.
\end{proof}

\begin{lem}
For $y,y' \in B$,
\[
  \bigl| d_B(y,y') - 
  d_{B_k} \bigl( \varphi_k (y), \varphi_k (y') \bigr) \Bigr| 
  \le C \sqrt{\frac {\log k}k} \,.
\]
\end{lem}

\begin{proof}
From the above discussion, it suffices to consider on $B(k)$.
Let $\gamma : [0,l] \to B(k)$ 
be a geodesic with $\gamma (0) = y$ and $\gamma (l) = y'$.
We take its horizontal lift
$\tilde{\gamma} : [0,l] \to X$ with 
$\tilde{\gamma} (0) = y \in B \subset X$ and 
$z := \tilde{\gamma} (l) \in \pi^{-1}(y')$.
From Lemma \ref{claim1},
\begin{align*}
  d_B(y,y') 
  &= \int_0^l \left| \frac{\dr}{\dr t} \tilde{\gamma} \right| dt \
   \ge \int_0^l 
   \left| \left( \frac{\dr}{\dr t} \tilde{\gamma} \right)^H \right| dt
   - \frac C{\sqrt k} \\
  & \ge d_{B_k}( \varphi_k (y) , \varphi_k (z)) - \frac C{\sqrt k}\\
  & \ge d_{B_k}( \varphi_k (y) , \varphi_k (y')) 
    - d_{B_k}( \varphi_k (y') , \varphi_k (z)) - \frac C{\sqrt k} \,.
\end{align*}
From the proof of Lemma \ref{claim3}, 
$|d_{B_k}( \varphi_k (y') , \varphi_k (z))| \le \frac C{\sqrt k}$.
Hence we obtain
\begin{equation}
  d_B(y,y') - d_{B_k}( \varphi_k (y) , \varphi_k (y')) 
  \ge \frac C{\sqrt k} \,.
  \label{estimate1}
\end{equation}

Next we take a geodesic $\gamma : [0,l] \to B_k(k)$ with
$\gamma (0) = \varphi_k(y)$ and $\gamma (l) = \varphi_k(y')$.
Let $\tilde{\gamma} : [0,l] \to X$ be a lift of $\gamma$.
Then 
\[
  d_{B_k}( \varphi_k (y) , \varphi_k (y'))
  = \int_0^l \left| \left( \frac{d}{d t} \tilde{\gamma} 
  \right)^H \right|_{\omega_k} dt \, .
\]
We write 
\[
  \frac{d}{dt} \tilde{\gamma} 
  = \xi + \eta,
  \quad \xi \in T_{X/B} \,, 
  \quad \eta \in (T_{X/B})^{\perp}\,.
\]
From Lemma \ref{claim1},
\[
  \left| \left( \frac{d}{dt} \tilde{\gamma} \right)^H
  - \eta \right| \le \frac C{\sqrt k} \,.
\]
Hence we have
\begin{align*}
  d_{B_k}( \varphi_k (y) , \varphi_k (y'))
  &\ge  \int_0^l |\eta| dt - \frac C{\sqrt k} \\
  &=  \int_0^l \left| \frac{d}{dt} \pi (\tilde{\gamma}) 
  \right|_{\omega} dt - \frac C{\sqrt k} \\
  &\ge d_{B}( \pi(\tilde{\gamma}(0)), \pi(\tilde{\gamma}(l))) 
  - \frac C{\sqrt k}\,.
\end{align*}
Since
\[
  d_{B}( \pi(\tilde{\gamma}(0)), y) \, , \,\,
  d_{B}( \pi(\tilde{\gamma}(l)), y') \le \frac C{\sqrt k}\,,
\]
we obtain
\begin{equation}
  d_{B_k}( \varphi_k (y) , \varphi_k (y'))
  \ge d_B(y,y') - \frac C{\sqrt k} \,.
  \label{estimate2}
\end{equation}
(\ref{estimate1}) and (\ref{estimate2}) prove the lemma.
\end{proof}

\subsection{Proof of Theorem \ref{main} (iii).}

From Theorem \ref{main} (ii), there exists a distance on 
$B \amalg B_k$ which is compatible with those on $B$ and $B_k$ and 
\[
  d(b ,\varphi_k(b)) \le C \sqrt{\frac {\log k}k}
\]
for any $b \in B$.

\begin{lem}\label{claim3-1}
For any $z \in X$,
\[
  d( \pi(z) , \pi_k(z) ) \le C \sqrt{\frac {\log k}k} \, .
\]
\end{lem}

\begin{proof}
From the proof of Lemma \ref{claim3}, 
\[
  d_{B_k}( \varphi_k (\pi(z)), \pi_k(z)) \le C \sqrt{\frac {\log k}k} \,.
\]
On the other hand, from the choice of the distance on $B \amalg B_k$,
\[
  d( \pi(x) ,\varphi_k(\pi(z)) ) \le C \sqrt{\frac {\log k}k} \, .
\]
Lemma \ref{claim3-1} follows from these two inequalities.
\end{proof}

For each $k$, we denote $ id_k := id_X : (X, \omega) \to (X, \omega_k)$.
Then, from the proof of Theorem \ref{GH_for_X}, 
there exists a distance on 
$(X,\omega) \amalg \left( \amalg_k (X, \omega_k) \right)$
which is compatible with $\omega$ and $\omega_k$ and
\[
  d(z, id_k(z)) \le C \sqrt{\frac {\log k}k}
\]
for any $z \in X$.

Take $z_k \in (X,\omega_k)$ converging to $z \in (X,\omega)$ 
in $(X,\omega) \amalg \left( \amalg_k (X, \omega_k) \right)$
as $k \to \infty$.
Then $id_k^{-1}(z_k) \to z$ in $(X,\omega)$
and we have
\[
  \pi \bigl( id_k^{-1}(z_k)\bigr) \longrightarrow \pi(z) \quad \text{in } B \,.
\]
Combining this with Lemma \ref{claim3-1}, we obtain
\[
  d(\pi_k(z_k), \pi(z)) \le 
  d(\pi_k(z_k), \pi(z_k)) + d_B(\pi(z_k), \pi(z)) \to 0
\]
as $k \to \infty$.
\qed



\begin{thebibliography}{99}

\bibitem{B}
  Th. Bauer,
  {\it Projective images of Kummer surfaces},
  Math. Ann. 299 (1994), 155--170.

\bibitem{DLM}
  X. Dai, K. Liu, and X. Ma,
  {\it On the asymptotic expansion of Bergman kernel},
  C. R. Math. Acad. Sci. Paris  339  (2004),  no. 3, 193--198.
  
\bibitem{F2}
  K. Fukaya,
  {\it Hausdorff convergence of Riemannian manifolds and its
   applications},
  Recent topics in differential and analytic geometry,  143--238,
  Adv. Stud. Pure Math., 18-{\rm I},
  Academic Press, Boston, MA, 1990. 
 
\bibitem{F}
  K. Fukaya,
  {\it Mirror symmetry of Abelian varieties and multi-theta
  functions},
  J. Alg. Geom. 11 (2002), 393--512.
  
\bibitem{M1}
  D. Mumford,
  {\it Tata lectures on theta. I},
  Progress in mathematics, 28, Birkhuser (1983).
  
\bibitem{N}
  Y. Nohara,
  {\it Projective embeddings and Lagrangian fibrations of
   abelian varieties},
  Math. Ann. 333 (2005), 741--757.

\bibitem{P}
  P. Petersen,
  {\it Riemannian geometry},
  Graduate Texts in Mathematics, 171. Springer-Verlag, New York, 1998.
  
\bibitem{PZ}
  A. Polishchuk and E. Zaslow,
  {\it Categorical mirror symmetry: the elliptic curve},
  Adv. Theor. Math. Phys. 2 (1998), 443--470.
  
\bibitem{Sa}
  R. Sasaki,
  {\it Bounds on the degree of the equations defining Kummer varieties},
  J. Math. Soc. Japan  33  (1981), no. 2, 323--333. 
  
\bibitem{S}
  J. Song,
  {\it The Szeg\"o kernel on an orbifold circle bundle},
  math.DG/0405071.
  
\bibitem{T}
  G. Tian,
  {\it On a set of polarized K\"ahler metrics on algebraic manifolds},
  J. Differential Geom. 32 (1990), 99--130.

\bibitem{W}
  J. Weitsman,
  {Quantization via real polarization of the moduli space of
  flat connections and Chern-Simons gauge theory in genus one},
  Comm. Math. Phys. 137 (1991), 175--190.
  
\bibitem{Zel}
  S. Zelditch,
  {\it Szeg\"o kernels and a theorem of Tian},
  International Math. Res. Notices 1998, No.6, 317--331.


\end{thebibliography}
\end{document}